\documentclass[a4paper,12pt]{article}
   \usepackage[top=2.5cm,bottom=2.5cm,left=2.5cm,right=2.5cm]{geometry}
    \usepackage{amssymb}
    \usepackage{graphicx}
    \usepackage{amsthm}
    \usepackage{amsmath}
    \usepackage[english]{babel}
    \usepackage[latin1]{inputenc}
    \usepackage{textcomp}
    \usepackage{longtable}

\newcommand{\floor}[1]{{\left\lfloor #1 \right\rfloor}}

\newtheorem{proposition}{Proposition}
\newtheorem{lemma}{Lemma}

\begin{document}

\title{\vspace*{2.5cm} Novel Bounds for the Normalized Laplacian Estrada and Normalized Energy Index of Graphs}
%\footnote{We thank Monica Bianchi and Anna Torriero for their valuable comments and suggestions.

\author{{Gian Paolo Clemente$^{a}$}, Alessandra Cornaro$^{a}$}
\date{}
\maketitle

\vspace*{-5mm}

\begin{center}
$^{{}^a}$\emph{Department of Mathematics and Econometrics, Catholic
University, Milan, Italy}\\[2mm]
{\tt
gianpaolo.clemente@unicatt.it, alessandra.cornaro@unicatt.it}\\[5mm]
\vspace{5mm}

\end{center}

\vspace{3mm}

\thispagestyle{empty}
\begin{abstract}
For a simple and connected graph, several lower and upper bounds of graph invariants expressed in terms of the eigenvalues of the normalized Laplacian matrix have been proposed in literature. In this paper, through a unified approach based on majorization techniques, we provide some novel inequalities depending on additional information on the localization of the eigenvalues of the normalized Laplacian matrix.
Some numerical examples show how sharper results can be obtained with respect to those existing in literature.
\end{abstract}
\baselineskip=0.30in

\newpage

\section{Introduction}\label{intro}
In literature, several topological indices, related to the structural properties of graphs, have been widely explored.
We focus here on the normalized Laplacian Estrada index (see \cite{Hakimi} and \cite{LiGuoShiu}) and the normalized Laplacian energy index (see \cite{Cavers}), that are based on a particular matrix associated with a graph, called the normalized Laplacian matrix. Properties about the spectrum of this matrix and its relationship to the Randi\'c index have been investigated in several works (see \cite{BroHae}, \cite{Cavers}, \cite{Chung} and \cite{Cve}).
In this paper we use a powerful methodology that relies on majorization techniques (see \cite{BCPT2}, \cite{BCPT3}, \cite{BCPT4} and \cite{BCT1}) in order to localize the graph topological indices we consider. In particular, through this technique, we derive new bounds for these indices taking advantage of additional information on the localization of the eigenvalues of normalized Laplacian matrix. Furthermore, this additional information can be quantified by using numerical approaches developed in \cite{CC1} and \cite{CC} and extended for normalized Laplacian matrix in \cite{CC2}.
Finally, some existing bounds (see \cite{Hakimi} and \cite{LiGuoShiu}), depending on well-known inequalities on Randi\'c index, have been also improved by using some novel results proposed in \cite{BCT2}.

\noindent The paper is organized as follows: in Section \ref{sec:Not} some preliminaries are given.
In Section 3 we provide, through majorization techniques, new bounds for topological indices expressed in terms of the eigenvalues of the normalized Laplacian matrix and we also recover in a straightforward way some results proposed in \cite{LiGuoShiu}. The relation between normalized Laplacian Estrada index and Randi\'c index has been used in Section \ref{sec:NEErand} to obtain new inequalities on normalized Laplacian Estrada index. Finally, in Section \ref{sec:NumRes} several numerical results are reported, showing how the proposed bounds are tighter than those given in literature.

\section{Notations and Preliminaries}\label{sec:Not}

\subsection{Basic graph concepts}\label{sec:gt}
\noindent We consider a simple, connected and undirected graph $G=(V,E)$ where $V=\{1, 2, \ldots, n\}$ is the set of vertices and $E\subseteq V\times V$ the set of edges, $|E|=m$.

\noindent The degree sequence of $G$ is denoted by $\pi=(d_{1},d_{2},..,d_{n})$ and it is arranged in non-increasing order $d_{1}\geq d_{2}\geq \cdots \geq d_{n}$, where $d_{i}$ is the degree of vertex $i$.

It is well known that $\overset{n}{\underset{i=1}{\sum }}d_{i}=2m$ and that if $G$ is a tree, i.e. a connected graph
without cycles, $m=n-1.$

Let $A(G)$ be the adjacency matrix of $G$ and $D(G)$ be the diagonal matrix of vertex degrees.
The matrix $L(G)=D(G)-A(G)$ is called Laplacian matrix of $G$, while $\mathcal{L}%
(G)=D(G)^{-1/2}L(G)D(G)^{-1/2}$ is known as normalized Laplacian.
Let $\lambda _{1}\geq \lambda _{2}\geq ...\geq \lambda _{n}$, $\mu _{1}\geq \mu _{2}\geq ...\geq \mu_{n}$ and $\gamma _{1}\geq \gamma_{2}\geq ...\geq \gamma_{n}$ be the set of (real) eigenvalues of $A(G)$, $L(G)$ and $\mathcal{L}(G)$ respectively.

We now recall some properties of normalized Laplacian eigenvalues useful for our purpose. For more details we refer the reader to \cite{BroHae}, \cite{Chung} and \cite{Cve}.

\begin{lemma}(see \cite{Chung})

Given a connected graph $G$ of order $n \geq 2$, the following properties of the spectrum of $\mathcal{L}(G)$ hold:

\begin{enumerate}

  \item $\overset{n}{\underset{i=1}{\sum }}\gamma_{i}=\mathrm{tr}(\mathcal{L}(G))=n;$
  \item $ \overset{n}{\underset{i=1}{\sum }}\gamma_{i}^{2}=\mathrm{tr}(\mathcal{L}^{2}(G))=n+2\sum_{(i,j)\in E}\frac{1}{d_{i}d_{j}};$
  \item $\frac{n}{n-1}\leq\gamma _{1}\leq 2.$ The left inequality is attained if and only if G is a complete graph, while the right inequality holds when G is a bipartite graph;
  \item $\gamma_{n}=0$, $\gamma_{n-1} \neq 0$ if $G$ is connected.
\end{enumerate}

\end{lemma}

\subsection{Normalized Laplacian indices}\label{subsec:tp}

The normalized Laplacian Estrada index has been proposed in \cite{LiGuoShiu} and it is defined as:
\begin{equation}
NEE(G)=\sum_{i=1}^{n}e^{\left(\gamma _{i}-1\right)}=\dfrac{1}{e}
\sum_{i=1}^{n}e^{\gamma _{i}}.
\label{eq:NEE}
\end{equation}

\noindent In \cite{Hakimi}, an alternative definition of normalized Laplacian Estrada index has been provided:
\begin{equation}
\ell EE(G)=\sum_{i=1}^{n}e^{\gamma _{i}}.
\label{eq:lEE}
\end{equation}
Notice that $NEE(G)=\frac{1}{e}\ell EE(G)$, any results derived for $NEE(G)$ can be trivially re-stated for $\ell EE(G)$ and viceversa.

Another graph invariant, introduced in \cite{Cavers}, is the normalized Laplacian energy index of a graph denoted by:
\begin{equation}
NE(G)=\sum_{i=1}^{n}\left|\gamma_{i}-1\right|.
\label{eq:En}
\end{equation}
%If $G$ is a $d-$regular graph, it is well-known that $E(G)=d \cdot NE(G)$.

\subsection{Randi\'c index and Majorization techniques}\label{subsec:mt}

The Randi\'c index is defined as:
\begin{equation*}
R_{-1}(G)=\underset{\left( i,j\right) \in E}{\sum }\left(
\dfrac{1}{d_{i}d_{j}}\right),
\end{equation*}

\noindent and it can be equivalently expressed as:
\begin{equation*}
R_{-1}(G)=\dfrac{1}{2}\left( \underset{\left(
i,j\right) \in E}{\sum }\left(\dfrac{1}{d_{i}}+\frac{1}{d_{j}}\right) ^{2}-\overset{n}{\underset{i=1}{\sum }}\dfrac{1}{d_{i}}\right).
\end{equation*}

\noindent Given a fixed degree sequence $\mathbf{\pi}$, let
$\mathbf{x}\in \mathbb{R}^{m}$ be the vector whose components are $\dfrac{1}{d_{i}}+\dfrac{1}{d_{j}},$
with $\left( i,j\right) \in E.$

%Since $\overset{n}{\underset{i=1}{\sum }}\dfrac{1}{d_{i}}$ is a constant, $R_{-1}(G)$ is a Schur convex
%function (for more details we refer the reader to \cite{Marshall}) of $\mathbf{x}$ and it is minimal (maximal) if and only if $f(\mathbf{x})=\overset{m}{\underset{i=1}{\sum }}x_{i}^{2}=$ $\left\Vert \mathbf{x
%}\right\Vert_{2}^{2}$ is minimal (maximal).
%It is easy to show that $\sum_{i=1}^{m}x_{i}=\underset{\left(i,j\right) \in E}{\sum }\left(\dfrac{1}{d_{i}}+\dfrac{1}{d_{j}}\right) =n$ and thus $\sum_{i=1}^{m}x_{i}$ is a constant.

\noindent Since $\sum_{i=1}^{m}x_{i}=\underset{\left(i,j\right) \in E}{\sum }\left(\dfrac{1}{d_{i}}+\dfrac{1}{d_{j}}\right) =n$,
let $\Sigma _{n}= \{\mathbf{x}\in \mathbb{R}_{+}^{m}: \sum_{i=1}^{m}x_{i} =n,  x_1 \ge x_2 \ge \cdots \ge x_m \}$.
By considering  a closed subset $S$ of $\Sigma_{n}$ whose maximal and minimal elements with respect to the majorization order are $\mathbf{x}^{\ast }(S)$ and $\mathbf{x}_{\ast }(S)$,  the Randi\'c index can be bounded as follows (see (5) in \cite{BCT2}):
\begin{equation}
L_{1}=\dfrac{\left\Vert \mathbf{x}_{\ast }(S)\right\Vert _{2}^{2}-\overset{n}{\underset{i=1}{
\sum }}\dfrac{1}{d_{i}}}{2} \leq R_{-1}(G)\leq \dfrac{\left\Vert \mathbf{x}^{\ast }(S)\right\Vert _{2}^{2}-\overset{n}{\underset{i=1}{%
\sum }}\dfrac{1}{d_{i}}}{2}=U_{1}.
\label{randicmaj}
\end{equation}

\noindent Inequalities (\ref{randicmaj}) will be used in Section \ref{sec:NEErand} in order to derive new bounds for $NEE(G)$.

Using the information available on the degree sequence of $G$ and characterizing the set $S$, the minimal and maximal elements $\mathbf{x}^{\ast}(S)$ and $\mathbf{x}_{\ast }(S)$ can be easily computed.

\noindent In this paper, we focus on a specific case of a graph $G$ with $h$ pendent vertices, whose degree sequence is of the type
\begin{equation}
\mathbf{\pi} =(d_{1},\cdots ,d_{n-h},\underbrace{1,\cdots ,1}_{h}),  \label{ds}
\end{equation}
where $h>0$ and $n-h\geq 2$ (we do not consider the star graph $S_{n}$ since it is well-known that $R_{-1}(S_{n})=1$).

\noindent It is noteworthy that this method could be applied to other suitable degree sequences.

Pointing out that $\dfrac{1}{d_{n-h}}+\dfrac{1}{d_{n-h-1}} < 1+\dfrac{1}{d_{1}}$ holds, we face the set

\begin{equation}
\begin{split}
S_{1}& =\left\{ \mathbf{x}\in \mathbb{R}_{+}^{m}:\sum_{i=1}^{m}x_{i} =n,\,\,\,1+\dfrac{1}{d_{1}}\leq x_{h}\leq \cdots
\leq x_{1}\leq \dfrac{1}{d_{n-h}}+1\right. , \\
& \left. \dfrac{1}{d_{1}}+\dfrac{1}{d_{2}}\leq x_{m}\leq \cdots \leq
x_{h+1}\leq \dfrac{1}{d_{n-h}}+\dfrac{1}{d_{n-h-1}}\right\}.
\end{split}
\label{S4n}
\end{equation}

For convenience of the reader, we report the expressions of the maximal and minimal elements of $S_{1}$.

%we report the main steps described in \cite{BCT2} and \cite{BCT1} to compute
%the maximal and minimal element with respect to the majorization order.

\noindent The maximal element is derived by means of Corollary 3  in \cite{BCT2} as follows:

\begin{equation}
\mathbf{x}^{\ast }(S_{1})=\left\{
\begin{array}{ccc}
\left[ \underset{k}{\underbrace{M_{1},.....,M_{1}}},\theta ,\underset{h-k-1}{
\underbrace{m_{1},.....,m_{1}},}\underset{m-h}{\underbrace{m_{2},.....,m_{2}}
}\right] & \text{ if } & n<a^{\ast } \\
&  &  \\
\left[ \underset{h}{\underbrace{M_{1},.....,M_{1}}},\underset{k-h}{
\underbrace{M_{2},.....,M_{2}}},\theta ,\underset{m-k-1}{\underbrace{
m_{2},.....m_{2}}}\right] & \text{ if } & n\geq a^{\ast }
\end{array},
\right.
\label{eq:maxel}
\end{equation}

\noindent where
\begin{equation*}
k=\left\{
\begin{array}{ccc}
\floor {\dfrac{n-h(m_{1}-m_{2})-mm_{2}}{M_{1}-m_{1}}} &
\text{ if } & n<a^{\ast } \\
&  &  \\
\floor {\dfrac{n-h(M_{1}-M_{2})-mm_{2}}{M_{2}-m_{2}}} &
\text{ if } & n\geq a^{\ast }
\end{array},
\right.
\end{equation*}

\noindent $a^{\ast }=hM_{1}+(m-h)m_{2}$, $m_{1}=1+\dfrac{1}{d_{1}}$, $m_{2}=\dfrac{1}{d_{1}}+\dfrac{1}{d_{2}}$, $M_{1}=1+\dfrac{1}{d_{n-h}}$, $M_{2}=\dfrac{1}{d_{n-h}}+\dfrac{1}{d_{n-h-1}}$ and $\theta$ is obtained as the difference between $n$ and the sum of the other components of the vector $\mathbf{x}^{\ast }(S_{1})$.

The minimal element is instead obtained by Corollary 10 in \cite{BCT2} as follows:
\begin{equation}
\mathbf{x}_{\ast }(S_{1})=\left\{
\begin{array}{ccc}
\left[\underset{h}{\underbrace{m_{1},...,m_{1}}},\underset{m-h}{\underbrace{\dfrac{n-hm_{1}}{m-h},...,\dfrac{n-hm_{1}}{m-h}}}\right] & \text{ if } & n<
\widetilde{a} \\
\left[\underset{h}{\underbrace{\dfrac{n-M_{2}(m-h)}{h},...,\dfrac{n-M_{2}(m-h)}{h}}},\underset{m-h}{\underbrace{M_{2},...,M_{2}}}\right] & \text{ if } &
n\geq \widetilde{a}, \\
&  &
\end{array},
\label{eq:minel}
\right.
\end{equation}
\noindent where $\widetilde{a}=hm_{1}+(m-h)M_{2}$ and $m_{1}$, $M_{2}$ have the same meaning of before.

\section{Bounds for normalized Laplacian indices via majorization techniques}

In this section we provide bounds for normalized Laplacian Estrada index and normalized Laplacian energy index.
These descriptors can be expressed in terms of Schur-convex or Schur-concave functions of suitable variables.
We briefly recall that Schur-convex (Schur-concave) functions preserve (reverse) the majorization order (see \cite{Marshall} for details).

\subsection{Normalized Laplacian Estrada index}\label{sec:NEEsub1}

Firstly, we focus on $NEE(G)$.
Let us consider the set
\begin{equation*}
S_{0}=\{ \mathbf{\gamma} \in \mathbb{R}^{n-1}: \overset{n-1}{\underset{i=1}{\sum }}\gamma_{i}= n,  \text{ } \gamma_{1} \geq \gamma_{2} \geq ... \geq \gamma_{n-2} \geq \gamma_{n-1} \geq 0\ \}.
\end{equation*}

\noindent We can now consider a subset $S^{1}_{0}$ of $S_{0}$:
\begin{equation*}
S^{1}_{0}=\{ \mathbf{\gamma} \in S_{0}: \gamma_{1} \geq \alpha\},
\end{equation*}
with $\alpha \geq \frac{n}{n-1}$.

In order to compute the minimal element of $S^{1}_{0}$, we apply Corollary 14 in \cite{BCT1} and we obtain:
\begin{equation*}
x_{\ast }(S^{1}_{0})=\left(\alpha, \underset{n-2}{\underbrace{\frac{n-\alpha}{n-2},...,\frac{n-\alpha}{n-2}}}\right).
\label{eq:minEl}
\end{equation*}

\noindent By the Schur-convexity of the function $NEE(G)$, we get the following bound:
\begin{equation}
%NEE(G)\geq \frac{1}{e}\left( 1+e^{\alpha}+(n-2)e^{\frac{n-alpha}{n-2}}\right) ,
NEE(G)\geq \frac{1}{e}+e^{\alpha-1}+(n-2)e^{\frac{2-\alpha}{n-2}}.
\label{eq:newB}
\end{equation}

\noindent Setting $\alpha = \dfrac{n}{n-1}$, we can easily derive the same result proved in \cite{LiGuoShiu}, Theorem 3.1:
\begin{equation}
NEE(G)\geq \left( n-1\right) e^{\dfrac{1}{n-1}}+\dfrac{1}{e}.
\label{eq:Libound}
\end{equation}

Furthermore, by applying a theoretical and numerical methodology (see \cite{BT} and \cite{CC2}), it is possible to compute a different lower bound $\alpha$ for the first eigenvalue of $\gamma_{1}$ in a fairly straightforward way, that is $\gamma_{1} \geq Q$, where
\begin{equation*}
Q=\frac{\left(n+\sqrt{\frac{b(h^{\ast}+1)-n^{2}}{h^{\ast}}}\right)}{(1+h^{\ast})},
\end{equation*}
with $b=n+2 \sum_{(i,j) \in E} \frac{1}{d_{i}d{j}}$ and $h^{\ast}=\left\lfloor {\frac{n^{2}}{b}}\right\rfloor$.

\noindent It is well-known that, for every connected graph of order $n$:
\begin{equation}
\left( \frac{2}{n}\right) \sum_{(i,j)\in E}\frac{1}{d_{i}d_{j}} \geq \frac{1}{n-1},
\label{eq:h}
\end{equation}
with inequality attained when $G \cong K_{n}$ (see \cite{BCPT1}). It has been shown in \cite{CC2} that $Q \geq \frac{n}{n-1}$ and thus
we assure that bound (\ref{eq:newB}), by placing $\alpha=Q$, is sharper than (\ref{eq:Libound}) (see \cite{BCT2} and \cite{BCT1} for more theoretical details).

We can further improve bound (\ref{eq:newB}) by identifying additional information on $\gamma_{2}$.
In this case we face the set:
\begin{equation*}
S^{2}_{0}=\{ \mathbf{\gamma} \in S_{0}: \gamma_{1} \geq \alpha\ \text{,} \gamma_{2} \geq \beta \}.
\end{equation*}
Under the assumptions $\alpha \geq \beta$ and $\alpha+\beta(n-2)>n$, by Corollary 14 in \cite{BCT1}, the minimal element of $S^{2}_{0}$ with respect to the majorization order is given by
\begin{equation*}
x_{\ast }(S^{2}_{0})=\left(\alpha,\beta, \underset{n-3}{\underbrace{\frac{n-\alpha-\beta}{n-3},...,\frac{n-\alpha-\beta}{n-3}}}\right)
\label{eq:minEl2}
\end{equation*}
and we can provide the following bound:
\begin{equation}
NEE(G)\geq \frac{1}{e}+e^{\alpha-1}+e^{\beta-1}+(n-3)e^{\frac{3-\alpha-\beta}{n-3}}.
\label{eq:newB2}
\end{equation}

In \cite{CC2} the authors found a lower bound $\beta$ for $\gamma_{2}$, that is $\gamma_{2}\geq R$, where: \\ $R=\frac{n-\sqrt{\frac{b(n-1)-n^{2}}{n-2}}}{n-1}.$
They proved that $R \leq Q$ and numerically showed, for some classes of graphs, that $Q+R(n-2)>n$, satisfying the conditions underlying Corollary 14 in \cite{BCT1}.

\noindent In virtue of these relations it is possible to compute bound (\ref{eq:newB2}) that is tighter than (\ref{eq:newB}) with $\alpha=Q$ and (\ref{eq:Libound}) (see \cite{BCT2} and \cite{BCT1} for more theoretical details).

Finally, for bipartite graphs, it is well-known that $\gamma_{1}=2$.
Hence
\begin{equation*}
S^{b}_{0}=\{ \mathbf{\gamma} \in \mathbb{R}^{n-2}: \overset{n-1}{\underset{i=2}{\sum }}\gamma_{i}= n-2, \text{ } 2 \geq \gamma_{2} \geq ... \geq \gamma_{n-2} \geq \gamma_{n-1} \geq 0\}.
\end{equation*}

\noindent By applying Corollary 14 in \cite{BCT1}, we recover the following bound provided in \cite{LiGuoShiu}, Theorem 3.2:
\begin{equation}
NEE(G)\geq \frac{1}{e}+ e+(n-2).
\label{Libound2}
\end{equation}

Also in this case, we can improve this bound by identifying additional information on $\gamma_{2}$.
We face the set:
\begin{equation*}
S^{2b}_{0}=\{ \mathbf{\gamma} \in S^{b}_{0}: \gamma_{2} \geq \beta \},
\end{equation*}
under the assumption $1 < \beta \leq 2$.
By Corollary 14 in \cite{BCT1}, the minimal element of $S^{2}_{0}$ with respect of majorization order is given by
\begin{equation*}
x_{\ast }(S^{2b}_{0})=\left(\beta, \underset{n-3}{\underbrace{\frac{n-2-\beta}{n-3},...,\frac{n-2-\beta}{n-3}}}\right)
\label{eq:minEl2b}
\end{equation*}
and we can provide the following bound:
\begin{equation}
NEE(G)\geq \frac{1}{e}+e+e^{\beta-1}+(n-3)e^{\frac{1-\beta}{n-3}},
\label{eq:newB2bip}
\end{equation}
where the lower bound $\beta=R$ of $\gamma_{2}$ derived in \cite{CC2} can be also used to compute (\ref{eq:newB2bip}).

In analogy with the results (\ref{eq:newB}) and (\ref{eq:newB2}) on $NEE(G)$, we can easily derive the following bounds for $\ell EE(G)$ for connected non bipartite graphs:
\begin{equation}
\ell EE(G)\geq 1+e^{\alpha}+(n-2)e^{\frac{n-\alpha}{n-2}}
\label{eq:newBl}
\end{equation}
and
\begin{equation}
\ell EE(G)\geq e^{\alpha}+e^{\beta}+(n-3)e^{\frac{n-\alpha-\beta}{n-3}},
\label{eq:newBl2}
\end{equation}
with $\gamma_{1} \geq \alpha$, $\gamma_{2} \geq \beta $ and $\alpha+\beta(n-2)>n$.

\noindent In Section \ref{sec:NumRes} we will compare these bounds with those proposed in \cite{Hakimi} and \cite{LiGuoShiu}.

\subsection{Normalized Laplacian energy index}\label{sec:NEG}

The normalized Laplacian energy index $NE(G)$ can be rewritten as a Schur-concave function of the variables $\left(\gamma_{i}-1\right)^{2}, i=1,\cdots ,n$:
\begin{equation}
NE(G)=1+\sum_{i=1}^{n-1}\sqrt{\left(\gamma_{i}-1\right)^{2}}.
\label{eq:En}
\end{equation}

\noindent If a lower bound for $\gamma_{1}$ is available, i.e. $%
\gamma_{1}\geq \alpha \left(\geq \dfrac{n}{n-1} \right)$, introducing the new variables $x_{i}=\left(\gamma_{i}-1\right)^{2}$ as a function of the eigenvalue $\gamma_{i}$ arranged in nonincreasing order, we get:
\begin{equation*}
x_{1}\geq k_{1} = \left(\alpha-1\right)^{2}.
\end{equation*}%
Let us consider the set
\begin{equation*}
S_{NE}= \{ \mathbf{x} \in \mathbb{R}^{n-1}: \overset{n-1}{\underset{i=1}{\sum }}x_{i}= 2\sum_{(i,j)\in E}\frac{1}{d_{i}d_{j}}-1, x_{1}\geq k_{1} \},
\label{setNE}
\end{equation*}%
where the relation $\overset{n-1}{\underset{i=1}{\sum }}x_{i}= 2\sum_{(i,j)\in E}\frac{1}{d_{i}d_{j}}-1$ has been obtained by using properties recalled in Lemma 1.

With the same methodology described for $NEE(G)$, we can derive the minimal element of $S_{NE}$ and then
the following upper bound:
\begin{equation}
NE(G)\leq 1+\sqrt{k_{1}}+\sqrt{(n-2)\left(a-k_{1}\right) },
\label{eq:NEQ}
\end{equation}
with $a=2\sum_{(i,j)\in E}\frac{1}{d_{i}d_{j}}-1$.
This bound could be computed by placing $k_{1}=\left(Q-1\right)^{2}$.

%By applying Corollary 14 in \cite{BCT1}, we get:
%\begin{equation*}
%\mathbf{x}_{\ast }(S_{NE})=\left[ k_{1},\underbrace{\dfrac{a-k_{1}}{n-2},\cdots ,%
%\dfrac{a-k_{1}}{n-2}}_{n-2}\right],
%\end{equation*}%

Considering also an additional information on $\gamma_{2}$ (i.e. $\gamma_{2}\geq \beta$), we may face the set:
\begin{equation*}
S^{2}_{NE}=\{ \mathbf{x} \in S_{NE}: x_{2} \geq k_{2} \}
\end{equation*}
under the assumptions $\alpha \geq \beta$ and $\alpha+\beta(n-2)>a$.

\noindent In this case, by means of the minimal element of $S^{2}_{NE}$, we can provide the bound:
\begin{equation}
NE(G)\leq 1+\sqrt{k_{1}}+\sqrt{k_{2}}+\sqrt{(n-3)\left( a-k_{1}-k_{2}\right) },
\label{eq:NEQR}
\end{equation}
where we can place $k_{1}=\left(Q-1\right)^{2}$ and $k_{2}=\left(R-1\right)^{2}$.

% is given by
%\begin{equation*}
%x_{\ast }(S^{2}_{NE})=\left(k_{1},k_{2}, \underset{n-3}{\underbrace{\frac{a-k_{1}-k_{2}}{n-3},...,\frac{a-k_{1}-k_{2}}{n-3}}}\right),
%\end{equation*}

Finally, for bipartite graphs, taking into account that $\gamma_{1}=2$, we set:
\begin{equation*}
S^{b}_{NE}=\{ \mathbf{x} \in \mathbb{R}^{n-2}: \overset{n-1}{\underset{i=2}{\sum }}x_{i}= 2\sum_{(i,j)\in E}\frac{1}{d_{i}d_{j}}-2\},
\label{setNE}
\end{equation*}%

\noindent and we derive the bound:
\begin{equation}
NE(G)\leq 2+\sqrt{a(n-2)},
\label{eq:NEQbip}
\end{equation}
where $a=2\sum_{(i,j)\in E}\frac{1}{d_{i}d_{j}}-2$.
\noindent
%By applying Corollary 14 in \cite{BCT1} where, in this case, $a=2\sum_{(i,j)\in E}\frac{1}{d_{i}d_{j}}-2$, we derive

Also in this case, we can improve this result by identifying additional information on $\gamma_{2}$.
We face the set:
\begin{equation*}
S^{2b}_{NE}=\{ \mathbf{x} \in S^{b}_{NE}: x_{2} \geq k_{2} \}
\end{equation*}
under the assumption $\dfrac{a-2}{n-2} < \beta \leq 2 $ and we can provide the bound:
\begin{equation}
NE(G)\leq 2+\sqrt{k_{2}}+\sqrt{(n-3)\left(a-k_{2}\right)},
\label{eq:NEQRbip}
\end{equation}
where the information $k_{2}=\left(R-1\right)^{2}$ can be used to compute (\ref{eq:NEQRbip}).

%By Corollary 14 in \cite{BCT1} the minimal element of $S^{2b}_{NE}$ with respect of majorization order is given by
%\begin{equation*}
%x_{\ast }(S^{2b}_{NE})=\left(k_{2}, \underset{n-3}{\underbrace{\frac{a-2-k_{2}}{n-3},...,\frac{a-2-k_{2}}{n-3}}}\right)
%\end{equation*}

\section{Bounds through Randi\'c Index}\label{sec:NEErand}

In Theorem 3.4 and Theorem 3.5 in \cite{LiGuoShiu}, the authors provided lower and upper bounds for $NEE(G)$ of a (bipartite) graph in terms of $n$ and maximum (or minimum) degree. This result has been obtained through well-known inequalities on Randi\'c index (see \cite{Shi}), i.e. $\dfrac{n}{2d_{1}} \leq R_{-1}(G) \leq \dfrac{n}{2d_{n}}.$

Following this idea, we now deduce some bounds for $NEE(G)$ and its variant $\ell EE(G)$ by using the methodology based on majorization recalled in Section \ref{subsec:mt}. In Section \ref{sec:NEErandres} we will numerically show that the bounds obtained are tighter than those provided in \cite{Hakimi} and \cite{LiGuoShiu}.

\noindent In virtue of (\ref{randicmaj}) and by means of Theorem 3.4 in \cite{LiGuoShiu}, we easily get the following result for bipartite graph:

\begin{proposition}
Let $G$ be a simple, connected and bipartite graph of order $n$. Then the normalized Laplacian Estrada index of $G$ is bounded as:

\begin{equation}
\frac{1}{e}+e+\sqrt{(n-2)^{2}+4(L_{1}-1)} \leq NEE(G) \leq \frac{1}{e}+e +(n-3)-\sqrt{2(U_{1}-1)}+e^{2(U_{1}-1)}.
\label{eq:newBrandBip}
\end{equation}

\end{proposition}

\noindent In the same way as before, by Theorem 3.5 in \cite{LiGuoShiu} we have the following bounds for non-bipartite graphs:
\begin{proposition}
Let $G$ be a simple and connected graph of order $n$. Then the normalized Laplacian Estrada index of $G$ is bounded as follows:

\begin{equation}
\sqrt{(n-1)(1+(n-2)e^\frac{2}{n-1})+4L_{1}} \leq NEE(G) \leq \frac{1}{e}+(n-1)-\sqrt{2U_{1}-1}+e^{2U_{1}-1}.
\label{eq:newBrand}
\end{equation}

\end{proposition}

\noindent Notice that, replacing $L_{1}=\dfrac{n}{2d_{1}}$ and  $U_{1}=\dfrac{n}{2d_{n}}$, we recover the same bounds provided in \cite{LiGuoShiu}, Theorem 3.4 and 3.5.

Bounds (\ref{eq:newBrandBip}) and (\ref{eq:newBrand}) can be trivially derived for $\ell EE(G)$ by using the proportionality relationship with $NEE(G)$.
For the comparisons provided in Section \ref{sec:NEErandres}, we only report the bound obtained for non-bipartite graph:
\begin{equation}
\sqrt{(n-1)(e+(n-2)e^\frac{n+1}{n-1})+4eL_{1}} \leq \ell EE(G) \leq 1+e\left[(n-1)-\sqrt{2U_{1}-1}\right]+e^{2U_{1}}.
\label{eq:EERand}
\end{equation}

\section{Numerical Results}\label{sec:NumRes}
\subsection{Comparing Bounds derived via majorization techniques}\label{sec:numMaj}
\subsubsection{Normalized Laplacian Estrada index}\label{sec:numMajNEE}
Firstly, we focus on $NEE(G)$ by comparing for non-bipartite graphs bounds (\ref{eq:newB}) and (\ref{eq:newB2}) with (\ref{eq:Libound}) proposed in \cite{LiGuoShiu}. It has been already analytically proved in Section \ref{sec:NEEsub1} that, when the additional information $\gamma_{1} \geq Q$ is considered, bound (\ref{eq:newB}) with $\alpha=Q$ is tighter than (\ref{eq:Libound}). We now show how these bounds behave according to different graphs.
In particular we analyze two alternative classes of graphs generated by using either the Erd\"os-R\'enyi (ER) model $G_{ER}(n,q)$ (see \cite{Boll}, \cite{Chung}, \cite{ER59} and \cite{ER60}) or the Watts and Strogatz (WS) model (see \cite{WS}).
Both models have been generated by using a well-known package of R (see \cite{igraph}) and by assuring that the graph obtained is connected.
The ER is constructed by connecting nodes randomly such that edges are included with probability $q$ independent from every other edge.
The WS networks have been derived beginning by a simulated $n$-node lattice and rewiring each edge at random to a new target node with probability $p$. As described by \cite{WS}, we choose a vertex and the edge that connects it to its nearest neighbor in a clockwise sense. With probability $p$, we reconnect this edge to a vertex chosen uniformly at random over the entire ring, with duplicate edges forbidden; otherwise
we leave the edge in place. We repeat this process by moving clockwise around the ring, considering each vertex in turn until one lap is completed. Next, we consider the edges that connect vertices to their second-nearest neighbors clockwise. As before, we randomly rewire each of these edges with probability $p$ and continue this process, circulating around the ring and proceeding outward to more distant neighbors after each lap, until each edge in the original lattice has been considered once. This construction allows to analyze the behavior of networks between regularity ($p = 0$) and disorder ($p = 1$).

In Table \ref{tab:res} we report the $NEE(G)$ index and the values of the three mentioned bounds evaluated on non-bipartite graphs generated by using ER model with different number of vertices and with $q$ equal to $0.5$.
Relative errors $r$ measures the absolute value of the difference between the lower bounds and $NEE(G)$ divided by the value of $NEE(G)$.

\begin{table}[!h]
\tiny
\centering
\begin{tabular}{|c||c||c|c|c||c|c|c|}
\hline\hline
$n$ & $NEE(G)$ & bound (\ref{eq:Libound}) & bound (\ref{eq:newB}) & bound (\ref{eq:newB2}) & r(\ref{eq:Libound})  & r(\ref{eq:newB})  & r(\ref{eq:newB2}) \\ \hline\hline

4	 &5.0862 	 &4.5547 	 &4.6783 	 &4.7112 	&10.4488\%	&8.0184\%	&7.3717\%\\ \hline
5	 &6.6073 	 &5.5040 	 &5.6407 	 &5.6935 	&16.6991\%	&14.6301\%	&13.8304\%\\ \hline
6	 &6.9783 	 &6.4749 	 &6.5088 	 &6.5265 	&7.2140\%	&6.7287\%	&6.4748\%\\ \hline
7	 &8.4965 	 &7.4560 	 &7.5345 	 &7.5559 	&12.2457\%	&11.3223\%	&11.0700\%\\ \hline
8	 &9.3463 	 &8.4428 	 &8.4778 	 &8.4933 	&9.6663\%	&9.2921\%	&9.1266\%\\ \hline
9	 &10.0295 	 &9.4331 	 &9.4456 	 &9.4541 	&5.9466\%	&5.8219\%	&5.7365\%\\ \hline
10	 &10.9027 	 &10.4256 	 &10.4334 	 &10.4391 	&4.3768\%	&4.3048\%	&4.2528\%\\ \hline
20	 &20.9252 	 &20.3947 	 &20.3963 	 &20.3977 	&2.5353\%	&2.5274\%	&2.5206\%\\ \hline
30	 &30.9411 	 &30.3853 	 &30.3860 	 &30.3867 	&1.7963\%	&1.7940\%	&1.7919\%\\ \hline
50	 &50.9236 	 &50.3782 	 &50.3784 	 &50.3786 	&1.0710\%	&1.0705\%	&1.0701\%\\ \hline
100	 &100.9001 	 &100.3729 	 &100.3730 	 &100.3731 	&0.5225\%	&0.5224\%	&0.5223\%\\ \hline \hline

\end{tabular}
\caption[]{Lower bounds for $NEE(G)$ and relative errors for graphs generated by $ER(n,0.5)$ model.}
\label{tab:res}
\end{table}

As expected, using bound (\ref{eq:newB2}) we observe an improvement with respect to existing bound according to all the analyzed graphs. The improvement is very significant for graphs with a small number of vertices, while it reduces for very large graphs. However, for large graphs formula (\ref{eq:Libound}) provided in \cite{LiGuoShiu} already gives a very low relative error.

The comparison has been extended in order to test the behaviour of the bounds on alternative graphs generated by using always the ER model with a different probability $q$. For sake of simplicity we report only the relative errors derived for graphs generated by using respectively $q=0.1$ and $q=0.9$ (see Table \ref{tab:res2}). In all cases bound (\ref{eq:newB2}) assures the best approximation to $NEE(G)$. We observe a best behaviour of all bounds when $q=0.9$ because we are moving towards the complete graph. We have indeed that the density of the graphs increases as long as greater probabilities are considered.

\begin{table}[!h]
\tiny
\centering
\begin{tabular}{|c||c|c|c|c||c|c|c|c|}
\hline\hline
 &\multicolumn{4}{c}{\textbf{q=0.1}} \vline  &\multicolumn{4}{c}{\textbf{q=0.9}}\\ \hline
$n$ & $NEE(G)$  & r(\ref{eq:Libound})  & r(\ref{eq:newB})  & r(\ref{eq:newB2}) & $NEE(G)$  & r(\ref{eq:Libound})  & r(\ref{eq:newB})  & r(\ref{eq:newB2})
\\ \hline\hline

4	 &5.3414 	&14.7282\%	&10.5879\%	&9.5649\%	 &4.5547 	&0	&0	&0\\ \hline
5	 &6.6073 	&16.6991\%	&14.6301\%	&13.8304\%	 &5.6685 	&2.9025\%	&2.5585\%	&2.0332\%\\ \hline
6	 &7.7763 	&16.7360\%	&14.3140\%	&13.9214\%	 &6.6355 	&2.4206\%	&2.2596\%	&2.1703\%\\ \hline
7	 &8.4997 	&12.2785\%	&11.7944\%	&11.5183\%	 &7.5599 	&1.3736\%	&1.3162\%	&1.27981\%\\ \hline
8	 &9.9938 	&15.5196\%	&14.3868\%	&14.1672\%	 &8.4664 	&0.2785\%	&0.2708\%	&0.2653\%\\ \hline
9	 &11.0383 	&14.5423\%	&13.8300\%	&13.6702\%	 &9.4702 	&0.3917\%	&0.3837\%	&0.3778\%\\ \hline
10	 &12.5449 	&16.8940\%	&15.8599\%	&15.7231\%	 &10.4541 	&0.2735\%	&0.2692\%	&0.2659\%\\ \hline
20	 &23.8531 	&14.4989\%	&14.1869\%	&14.1557\%	 &20.4488 	&0.2645\%	&0.2637\%	&0.2630\%\\ \hline
30	 &34.8955 	&12.9249\%	&12.7815\%	&12.7693\%	 &30.4424 	&0.1876\%	&0.1874\%	&0.1872\%\\ \hline
50	 &54.7998 	&8.0688\%	&8.0399\%	&8.0370\%	 &50.4274 	&0.0977\%	&0.0976\%	&0.0976\%\\ \hline
100	 &104.7347 	&4.1646\%	&4.1608\%	&4.1604\%	 &100.4293  &0.0562\%	&0.0561\%	&0.0561\%\\ \hline \hline

\end{tabular}
\caption[]{Lower bounds for $NEE(G)$ and relative errors for graphs generated respectively by $ER(n,0.1)$ and $ER(n,0.9)$ models.}
\label{tab:res2}
\end{table}

Finally, graphs have been simulated by using WS model with different rewiring probabilities $p$. As well-known, intermediate values of $p$ result in small-world networks that share properties of both regular and random graphs. In \cite{WS}, the authors show that these networks have small mean path lengths and high clustering coefficients. There is indeed a broad interval of $p$ over which the average path is almost as small as random yet the clustering coefficient is significantly greater than random. These small-world networks result from the immediate drop in average path caused by the introduction of few long-range edges.
In particular, we analyze the behaviour of bounds in this interval by considering graphs generated with a rewiring probability in the range $p \in (0.01,0.1)$. At this regard, Table \ref{tab:res3} reports bounds evaluated by considering $p=0.1$. In this case, we observe greater relative errors especially for large graphs. Probably, being these networks very far from complete graphs, bounds tend to assure a weaker approximation.
Similar results have been obtained by simulating WS graphs choosing different values of $p$ that belong to the interval.

\begin{table}[!h]
\tiny
\centering
\begin{tabular}{|c||c||c|c|c||c|c|c|}
\hline\hline
$n$ & $NEE(G)$ & bound (\ref{eq:Libound}) & bound (\ref{eq:newB}) & bound (\ref{eq:newB2}) & r(\ref{eq:Libound})  & r(\ref{eq:newB})  & r(\ref{eq:newB2}) \\ \hline\hline

4	 &5.0862 	 &4.5547 	 &4.6783 	 &4.7112 	&10.4488\%	&8.0184\%	&7.3717\%\\ \hline
5	 &6.3276 	 &5.5040 	 &5.6002 	 &5.6389 	&13.0165\%	&11.4961\%	&10.8843\%\\ \hline
6	 &7.5967 	 &6.4749 	 &6.5492 	 &6.5856 	&14.7666\%	&13.7886\%	&13.3087\%\\ \hline
7	 &8.8273 	 &7.4560 	 &7.6023 	 &7.6273 	&15.5347\%	&13.8778\%	&13.5948\%\\ \hline
8	 &10.1431 	 &8.4428 	 &8.5646 	 &8.5878 	&16.7630\%	&15.5621\%	&15.3339\%\\ \hline
9	 &11.3946 	 &9.4331 	 &9.5349 	 &9.5568 	&17.2145\%	&16.3209\%	&16.1287\%\\ \hline
10	 &12.6329 	 &10.4256 	 &10.5736 	 &10.5917 	&17.4727\%	&16.3005\%	&16.1573\%\\ \hline
20	 &25.2327 	 &20.3947 	 &20.5345 	 &20.5442 	&19.1737\%	&18.6195\%	&18.5813\%\\ \hline
30	 &37.7967 	 &30.3853 	 &30.5175 	 &30.5240 	&19.6086\%	&19.2590\%	&19.2418\%\\ \hline
50	 &63.2448 	 &50.3782 	 &50.5227 	 &50.5268 	&20.3442\%	&20.1157\%	&20.1092\%\\ \hline
100	 &126.4764 	 &100.3729 	 &100.5201 	 &100.5222 	&20.6390\%	&20.5226\%	&20.5209\%\\ \hline\hline

\end{tabular}
\caption[]{Lower bounds for $NEE(G)$ and relative errors for graphs generated by $WS(n,0.1)$ model.}
\label{tab:res3}
\end{table}

\subsubsection{Normalized Laplacian energy index}\label{sec:numMajNE}
We compare here bounds proposed in Section \ref{sec:NEG} for $NE(G)$ with the following upper bounds proposed in \cite{Cavers}:
\begin{equation}
NE(G) \leq 2  \left\lfloor\frac{n}{2}\right\rfloor,
\label{eq:Cav1}
\end{equation}
\begin{equation}
NE(G) \leq \sqrt{\frac{15}{28}} (n+1).
\label{eq:Cav2}
\end{equation}

Table \ref{tab:res7} reports main results derived for graphs generated by a $ER(n,0.5)$ model.
We observe how both bounds (\ref{eq:NEQ}) and (\ref{eq:NEQR}) are tighter than those proposed in \cite{Cavers}.
The improvement increases for greater number of vertices.

\begin{table}[!h]
\tiny
\centering
\begin{tabular}{|c||c||c|c||c|c||}
\hline\hline
$n$ & $NE(G)$ & bound (\ref{eq:Cav1}) & bound (\ref{eq:Cav2}) & bound (\ref{eq:NEQ})& bound (\ref{eq:NEQR})  \\ \hline\hline

4	 &3.00 	&4	 &3.66 	 &3.12 	 &3.05\\ \hline
5	 &2.55 	&4	 &4.39 	 &2.70 	 &2.61\\ \hline
6	 &3.15 	&6	 &5.12 	 &3.87 	 &3.62\\ \hline
7	 &3.81 	&6	 &5.86 	 &4.51 	 &4.27\\ \hline
8	 &4.32 	&8	 &6.59 	 &4.69 	 &4.47\\ \hline
9	 &3.90 	&8	 &7.32 	 &4.34 	 &4.14\\ \hline
10	 &3.58 	&10	 &8.05 	 &4.00 	 &3.83\\ \hline
20	 &5.01 	&20	 &15.37 	 &5.68 	 &5.56\\ \hline
30	 &5.60 	&30	 &22.69 	 &6.43 	 &6.33\\ \hline
50	 &7.31 	&50	 &37.33 	 &8.44 	 &8.36\\ \hline
100	 &9.59 	&100 &73.92 	 &11.13  &11.08\\ \hline\hline

\end{tabular}
\caption[]{Upper bounds for $NE(G)$ for graphs generated by $ER(n,0.5)$ model.}
\label{tab:res7}
\end{table}

Considering instead WS networks, derived as in Section \ref{sec:numMajNEE} by assuming a rewiring probability equal to $0.1$, we observe in Table \ref{tab:res8} greater values of $NE(G)$. In this case, bound (\ref{eq:Cav2}) gives better results than those observed for ER graphs. However it is confirmed the best approximation when bound (\ref{eq:NEQR}) is used.

\begin{table}[!h]
\tiny
\centering
\begin{tabular}{|c||c||c|c||c|c||}
\hline\hline
$n$ & $NE(G)$ & bound (\ref{eq:Cav1}) & bound (\ref{eq:Cav2}) & bound (\ref{eq:NEQ})& bound (\ref{eq:NEQR})  \\ \hline\hline
4	 &2 	&4	 &3.66 	 &2.72 	 &2.41\\ \hline
5	 &3.24 	&4	 &4.39 	 &3.45 	 &3.37\\ \hline
6	 &4 	&6	 &5.12 	 &4.15 	 &4.08\\ \hline
7	 &4.49 	&6	 &5.86 	 &4.87 	 &4.63\\ \hline
8	 &5.12 	&8	 &6.59 	 &5.57 	 &5.33\\ \hline
9	 &5.76 	&8	 &7.32 	 &6.27 	 &6.03\\ \hline
10	 &6.47 	&10	 &8.05 	 &6.99 	 &6.75\\ \hline
20	 &11.97 	&20	 &15.37 	 &13.67 	 &13.42\\ \hline
30	 &19.24 	&30	 &22.69 	 &21.16 	 &20.90\\ \hline
50	 &31.79 	&50	 &37.33 	 &35.27 	 &34.99\\ \hline
100	 &63.21 	&100	& 73.92 	& 70.22 	& 69.93\\ \hline \hline
\end{tabular}
\caption[]{Upper bounds for $NE(G)$ for graphs generated by $WS(n,0.1)$ model.}
\label{tab:res8}
\end{table}

\subsection{Bounds based on Randi\'c Index}\label{sec:NEErandres}

We now consider an example based on a specific degree sequence of type (\ref{ds}) in order to explain the details of the procedure used to bound $NEE(G)$ via Randi\'c Index. In the next we will extend the results to several degree sequences of type (\ref{ds}).

\noindent \textbf{Example 1}. Let us consider the class $C_{\pi}$ of graphs with the following degree sequence: \\
$\pi =  (7,6,5,4,4,4,3,3,3,3,3,3,2,2,2,2,1,1,1,1)$

%\begin{figure}[!hb]
%\centering
%\includegraphics[height=6cm]{GraficoRandic.jpg}
%\caption{Generated Graph}
%\label{fig:1}
%\end{figure}

\noindent We have $n = 20$, $m = 30$ and $h=4$ pendant nodes.
Since $\tilde{a} > n$, the minimal element (\ref{eq:minel}) is:

\begin{equation*}
\mathbf{x}_{\ast }(G)=\left[ \underset{4}{\underbrace{\frac{8}{7},...,
\frac{8}{7}}},\underset{16}{\underbrace{\frac{54}{91},...,\frac{54}{91}}}\right].
\end{equation*}

\noindent Replacing these values in (\ref{randicmaj}), we find $L_{1}=2.56$, while $\dfrac{n}{2d_{1}}=1.43$.

The bounds for $NEE(G)$ are figured out in Table 6. Furthermore, in order to test how these bounds behave, the exact value of $NEE(G)$ is also needed. Having a huge number\footnote{We estimate the total number of graphs with the degree sequence $\pi$ by using the importance sampling algorithm proposed in \cite{BlitzDiac}. Authors show robust results by applying the algorithm with 100.000 trials. In this case we derive a total number of graphs equal to roughly $1.20 \cdot 10^{20}$ with a standard error of $4 \cdot 10^{17}$. However it is noteworthy that also graphs belonging to the same isomorphism class are considered in this value. For the computation of the average values in Tables 6 and 7, we take into account only graphs with a different $NEE(G)$.} of graphs $G \in C_{\pi}$, we randomly generate one million of different graphs belonging to the class $C_{\pi}$. The average value, the minimum and maximum values of the index are also reported in Table 6.
\begin{table}[h]
\centering
\begin{tabular}{|l|l|}
\hline
Reference & Bound \\ \hline
Theorem 3.5 of \cite{LiGuoShiu} & $20.12$ \\ \hline
(\ref{eq:newBrand}) & $20.23$ \\ \hline
Min($NEE(G)$) & $20.51$ \\ \hline
Mean($NEE(G)$) & $23.25$ \\ \hline
Max($NEE(G)$) & $25.52$ \\ \hline
\end{tabular}
\label{tab:Tablelow}
\caption{Lower bounds for $NEE(G)$.}
\end{table}

%By using the first left inequality of (\ref{eq:newBrand}), we have $NEE(G) \geq 20.23$ against an exact value of $23.06$.
%Formula proposed in Theorem 3.5 of \cite{LiGuoShiu} leads to a lower bound of $20.12$.

Considering instead the upper bound, since $a^{*}< n$, we compute $k=12 > h=4$ and
we have that the maximal element (\ref{eq:maxel}) is:

\begin{equation*}
\mathbf{x}^{\ast }(G)=\left[ \underset{4}{\underbrace{\frac{3}{2},...,
\frac{3}{2}}},\underset{8}{\underbrace{1,...,
1}},\frac{31}{42},\underset{17}{\underbrace{\frac{31}{42},...,\frac{31}{42}}}\right],
\end{equation*}

\noindent leading to $U_{1}=4.96$, while $\dfrac{n}{2d_{n}}=10$.

Upper bounds and values of $NEE(G)$ are summarized in Table 7.

\begin{table}[h]
\centering
\begin{tabular}{|l|l|}
\hline
Reference & Bound \\ \hline
Theorem 3.5 of \cite{LiGuoShiu} & $1.7 * 10^{8}$ \\ \hline
(\ref{eq:newBrand}) & $7541.32$ \\ \hline
Min($NEE(G)$) & $20.51$ \\ \hline
Mean($NEE(G)$) & $23.25$ \\ \hline
Max($NEE(G)$) & $25.52$ \\ \hline
\end{tabular}
\label{tab:Tableup}
\caption{Upper bounds for $NEE(G).$}
\end{table}

\noindent Finally, if we know the value of Randi\'c Index, we can directly use it to compute (\ref{eq:newBrand}). For example, considering a random graph $G \in C_{\pi}$, we obtain  $R_{-1}(G)=3.0376$ deriving a better approximation (i.e. $20.27 \leq NEE(G) \leq 177.15$).

We now evaluate these bounds by randomly generating several degree sequences of type (\ref{ds}). For this aim, $ER(n,p)$ model has been used to derive different random graphs, where we disregard graphs whose degree sequence does not belong to the set (\ref{S4n}). The number of pendant vertices $h$ varies according to the specific degree sequence obtained.
Results have been compared to those analyzed in previous Section \ref{sec:numMajNEE}.
In particular, we report in Table \ref{tab:res4} bound (\ref{eq:Libound}) and bound (\ref{eq:newBrand}$d_{1}$) proposed in \cite{LiGuoShiu}, where bound (\ref{eq:newBrand}$d_{1}$) has been derived by using in (\ref{eq:newBrand}) the lower bound $\dfrac{n}{2d_{1}}$ of $R_{-1}(G)$. These bounds have been compared with bound (\ref{eq:newB}) and bound (\ref{eq:newB2}) already analysed in previous section and with bound (\ref{eq:newBrand}$L_{1}$) and bound (\ref{eq:newBrand}$R_{-1}$) evaluated by using the first left inequality of (\ref{eq:newBrand}) and by considering respectively the value of $L_{1}$ or by assuming to know the value of Randi\'c Index $R_{-1}(G)$.

We further observe that bound (\ref{eq:newB2}) based on value of $Q$ and $R$ shows the tighter lower bound in all cases by allowing a best approximation respect to bounds based on inequality (\ref{eq:newBrand}).
Furthermore, when inequality (\ref{eq:newBrand}) is considered, $L_{1}$ leads to a better bound than $\dfrac{n}{2d_{1}}$ used in \cite{LiGuoShiu}.
Finally, considering the exact value of Randi\'c Index we only get a slight improvement.

\begin{table}[!h]
\tiny
\centering
\begin{tabular}{|c|c|c||c||c|c||c|c|c|c|}
\hline\hline
$n$ & $m$ & $d_{1}$ & $NEE(G)$ & bound (\ref{eq:Libound}) & bound (\ref{eq:newBrand}$d_{1}$) & bound (\ref{eq:newB}) & bound (\ref{eq:newB2}) & bound (\ref{eq:newBrand}$L_{1}$) & bound (\ref{eq:newBrand}$R_{-1}$)\\ \hline\hline

4	&4	&3	 &4.8846 	 &4.5547 	 &4.1657 	 &4.6466 	 &4.6718 	 &4.2797 	 &4.2840\\ \hline
5	&5	&3	 &6.2381 	 &5.5040 	 &5.2075 	 &5.6002 	 &5.6389 	 &5.3288 	 &5.3651\\ \hline
6	&9	&5	 &6.8224 	 &6.4749 	 &6.1022 	 &6.4977 	 &6.5099 	 &6.1841 	 &6.1952\\ \hline
7	&11	&5	 &7.9613 	 &7.4560 	 &7.1182 	 &7.4774 	 &7.4901 	 &7.1933 	 &7.2062\\ \hline
8	&12	&6	 &9.2439 	 &8.4428 	 &8.0967 	 &8.4667 	 &8.4817 	 &8.2018 	 &8.2463\\ \hline
9	&14	&4	 &10.4376 	 &9.4331 	 &9.1872 	 &9.4722 	 &9.4843 	 &9.2432 	 &9.2490\\ \hline
10	&15	&8	 &11.2359 	 &10.4256 	 &10.0706 	 &10.4435 	 &10.4522 	 &10.1687 	 &10.1977\\ \hline
20	&30	&7	 &23.3079 	 &20.3947 	 &20.1166 	 &20.4401 	 &20.4465 	 &20.2346 	 &20.2914\\ \hline
30	&45	&7	 &34.8334 	 &30.3853 	 &30.1255 	 &30.4342 	 &30.4385 	 &30.2355 	 &30.2914\\ \hline
50	&75	&6	 &58.8047 	 &50.3782 	 &50.1563 	 &50.4387 	 &50.4416 	 &50.2535 	 &50.3259\\ \hline
100	&150 &8	 &117.1121 	 &100.3729 	 &100.1199 	 &100.4310 	 &100.4325 	 &100.2397 	 &100.3214\\ \hline\hline
\end{tabular}
\caption[]{Lower bounds for $NEE(G)$.}
\label{tab:res4}
\end{table}

On the same graphs upper bounds have been also evaluated by using the right part of inequality (\ref{eq:newBrand}). We observe in Table \ref{tab:res5}
a huge approximation, especially for large graphs, when we apply formula proposed in \cite{LiGuoShiu} based on the upper bound $\dfrac{n}{2d_{n}}$ of Randi\'c Index $R_{-1}(G)$ (see bound (\ref{eq:newBrand}$d_{n}$)).
By considering the upper bound based on $U_{1}$ we are able to improve the results, but for large graphs we derive useless bounds in this case too.
We have indeed that even when we directly use the value of $R_{-1}(G)$ we derive bounds significantly larger for graphs with a great number of vertices.

\begin{table}[!h]
\tiny
\centering
\begin{tabular}{|c|c|c||c||c|c|c|}
\hline\hline
$n$ & $m$ & $d_{1}$ & $NEE(G)$ & bound (\ref{eq:newBrand}$d_{n}$) & bound (\ref{eq:newBrand}$U_{1}$) & bound (\ref{eq:newBrand}$R_{-1}$)\\ \hline\hline

4	&4	&3	 &4.8846 	 & 21.72 	 &5.09 	     &4.86\\ \hline
5	&5	&3	 &6.2381 	 & 56.97 	 &7.62 	     &7.62\\ \hline
6	&9	&5	 &6.8224 	 &151.55 	 &8.15 	     &6.65\\ \hline
7	&11	&5	 &7.9613 	 & 407.35 	 &10.32 	 &8.16\\ \hline
8	&12	&6	 &9.2439 	 &1,101.36 	 &16.72 	 &10.86\\ \hline
9	&14	&4	 &10.4376 	 &2,986.50 	 &14.15 	 &13.19\\ \hline
10	&15	&8	 &11.2359 	 &8,109.45 	 &43.95 	 &12.78\\ \hline
20	&30	&7	 &23.3079 	 &1.78E+08	 &8,371.34 	 &236.01\\ \hline
30	&45	&7	 &34.8334 	 &3.93E+12	 &1.42E+06	&4.04E+03\\ \hline
50	&75	&6	 &58.8047 	 &1.91E+21	 &5.44E+09	&7.70E+06\\ \hline
100	&150 &8	 &117.1121 	 &9.89E+42	 &2.25E+22	&5.79E+13\\ \hline\hline
\end{tabular}
\caption[]{Upper bounds for $NEE(G)$.}
\label{tab:res5}
\end{table}

Bounds proposed for $\ell EE(G)$ have been also compared to the following bounds presented in \cite{Hakimi}:
\begin{equation}
\ell EE(G)> n e,
\label{eq:Hak1}
\end{equation}
\begin{equation}
\ell EE(G)> 2+\sqrt{n(n-1)e^{2}-6n+4},
\label{eq:Hak2}
\end{equation}
\begin{equation}
\ell EE(G)> \sqrt{n(n-1)e^{2}+4R_{-1}(G)+5n}.
\label{eq:Hak3}
\end{equation}

We observe in Table \ref{tab:res6} how the proposed bounds significantly improve those in \cite{Hakimi}.

\begin{table}[!h]
\tiny
\centering
\begin{tabular}{|c|c|c||c||c|c|c||c|c|c|}
\hline\hline
$n$ & $m$ & $d_{1}$ & $\ell EE(G)$ & bound (\ref{eq:Hak1}) & bound (\ref{eq:Hak2}) & bound (\ref{eq:Hak3})& bound (\ref{eq:newBl}) & bound (\ref{eq:newBl2}) & bound (\ref{eq:EERand}) \\ \hline\hline

4	&4	&3	&13.278	 &10.873 	 &6.173 	 &10.599 	 &12.631 	 &12.699 	 &11.633\\ \hline
5	&5	&3	&16.957	 &13.591 	 &7.636 	 &13.333 	 &15.223 	 &15.328 	 &14.485\\ \hline
6	&9	&5	&18.545	 &16.310 	 &7.991 	 &15.975 	 &17.663 	 &17.696 	 &16.810\\ \hline
7	&11	&5	&21.641	 &19.028 	 &9.265 	 &18.692 	 &20.326 	 &20.360 	 &19.553\\ \hline
8	&12	&6	&25.128	 &21.746 	 &10.173 	 &21.422 	 &23.015 	 &23.056 	 &22.295\\ \hline
9	&14	&4	&28.372	 &24.465 	 &12.482 	 &24.138 	 &25.748 	 &25.781 	 &25.126\\ \hline
10	&15	&8	&30.542	 &27.183 	 &12.059 	 &26.834 	 &28.388 	 &28.412 	 &27.642\\ \hline
20	&30	&7	&63.358	 &54.366 	 &24.421 	 &54.043 	 &55.562 	 &55.579 	 &55.003\\ \hline
30	&45	&7	&94.687	 &81.548 	 &36.926 	 &81.222 	 &82.729 	 &82.740 	 &82.189\\ \hline
50	&75	&6	&159.848	 &135.914 	 &66.935 	 &135.598 	 &137.107 	 &137.114 	 &136.603\\ \hline
100	&150 &8	&318.344	 &271.828 	 &133.313 	 &271.509 	 &273.000 	 &273.004 	 &272.480\\ \hline\hline

\end{tabular}
\caption[]{Lower bounds for $\ell EE(G)$.}
\label{tab:res6}
\end{table}

Considering instead the upper bounds, we compare our results with the following one in \cite{Hakimi}:
\begin{equation}
\ell EE(G)< e^{n} + R_{-1}(G) +\frac{n}{2}(3-n)-1.
\label{eq:Hak4}
\end{equation}

As reported in Table \ref{tab:res9}, upper bound (\ref{eq:EERand}) allows a better approximation than (\ref{eq:Hak4}).
Also in this case, the upper bounds do not show a good behaviour for large graphs.

\begin{table}[!h]
\tiny
\centering
\begin{tabular}{|c|c|c||c||c|c|}
\hline\hline
$n$ & $m$ & $d_{1}$ & $\ell EE(G)$ & bound (\ref{eq:Hak4}) & bound (\ref{eq:EERand}) \\ \hline\hline
4	&4	&3	 &13.28 	& 52.51 	& 13.83\\ \hline
5	&5	&3	 &16.96 	& 143.66 	& 20.73\\ \hline
6	&9	&5	 &18.55 	& 394.31 	& 22.15\\ \hline
7	&11	&5	 &21.64 	& 1,082.65 	& 28.06\\ \hline
8	&12	&6	 &25.13 	& 2,961.24 	& 45.46\\ \hline
9	&14	&4	 &28.37 	&8.08E+03	& 38.46\\ \hline
10	&15	&8	 &30.54 	&2.20E+04	& 119.46\\ \hline
20	&30	&7	 &63.36 	&4.85E+08	&2.28E+04\\ \hline
30	&45	&7	 &94.69 	&1.07E+13	&3.87E+06\\ \hline
50	&75	&6	 &159.85 	&5.18E+21	&1.48E+10\\ \hline
100	&150 &8	 &318.34 	&2.69E+43	&6.11E+22\\ \hline  \hline
\end{tabular}
\caption[]{Upper bounds for $\ell EE(G)$.}
\label{tab:res9}
\end{table}

\section{Conclusions}\label{Conc}
By using an approach for localizing some relevant graph topological indices based on the optimization of Schur-convex or Schur-concave functions, we derive some new bounds for normalized Laplacian Estrada index and for normalized Laplacian energy index.
The proposed bounds can be computed by using additional information on the localization of first and second eigenvalue of normalized Laplacian matrix.
A numerical section shows how this approach allows to derive tighter bounds than those provided in the literature. In particular, bound derived directly via majorization technique appear sharper than those depending by the Randi\'c Index. According to the latter ones, it is noteworthy that we analyzed only the results for a specific type of degree sequence, while different bounds could be derived for other suitable degree sequences.

\section*{Acknowledgement}
The authors are grateful to Monica Bianchi and Anna Torriero for useful advice and suggestions.

\section*{Competing interests}
The authors declare that they have no competing interests.

\section*{Authors' contribution}
All authors contributed equally to the writing of this paper. All authors read and approved the final manuscript.

\bibliographystyle{plain}
\bibliography{biblio}

\end{document}